\documentclass{ws-procs9x6}
\usepackage{amsmath}
\usepackage{amsfonts}

\newcommand{\tr}{{\rm tr}}

\newfont{\got}{eufm9 scaled 1095}

\newfont{\w}{msbm9 scaled\magstep1}
\def\R{\mbox{\w R}}

\newcommand{\const}{\rm{const}}

\begin{document}

\title{ON THE CURVATURE PROPERTIES OF REAL TIME-LIKE HYPERSURFACES OF K%
\"{A}HLER MANIFOLDS WITH NORDEN METRIC}

\author{M. MANEV, M. TEOFILOVA}

\address{Department of Geometry, Faculty of Mathematics and Informatics,\\ University of Plovdiv,
236 Bulgaria Blvd., Plovdiv 4003, Bulgaria,\\ e-mail:
mmanev@uni-plovdiv.bg, marta.teofilova@gmail.com}

\begin{abstract}
A type of almost contact hypersurfaces with Norden metric of a
K\"{a}hler manifold with Norden metric is considered. The curvature
tensor and the special sectional curvatures are characterized. The
canonical connection on such manifolds is studied and the form of
the corresponding K\"{a}hler curvature tensor is obtained. Some
curvature properties of the manifolds belonging to the widest
integrable main class of the considered type of hypersurfaces are
given.\footnote{ 
Hypersurface, curvature, almost contact manifold, Norden metric.}
\end{abstract}
\bodymatter

\section*{Introduction}

The K\"{a}hler manifolds with Norden metric have been introduced in
\cite {nor:A-spa}. These manifolds form the special class
$\mathcal{W}_{0}$ in the decomposition of the almost complex
manifolds with Norden metric, given in \cite {gabo:note}. This most
important class is contained in each of the basic classes in the
mentioned classification.

The natural analogue of the almost complex manifolds with Norden
metric in the odd dimensional case are the almost contact manifolds
with Norden metric, classified in \cite{gamigri:alcont}.

In \cite{man:hyp} two types of hypersurfaces of an almost complex
manifold with Norden metric are constructed as almost contact
manifolds with Norden metric, and the class of these hypersurfaces
of a $\mathcal{W}_{0}$-manifold is determined.

An important problem in the differential geometry of the K\"{a}hler
manifolds with Norden metric is the studying of the manifolds of
constant totally real sectional curvatures \cite{gagrimi:holohyp}.

In this paper we study some curvature properties of the real
time-like hypersurfaces of K\"{a}hler manifolds with Norden metric
of constant totally real sectional curvatures and particularly
curvature properties of nondegenerate special sections.

\section{Preliminaries}

\subsection{Almost complex manifolds with Norden metric}

Let $(M^{\prime},J,g^{\prime})$ be a $2n^{\prime}$-dimensional
almost complex manifold with Norden metric, i.e. $J$ is an almost
complex structure and $g^{\prime}$ is a metric on $M^{\prime}$ such
that:
\begin{equation*}
J^{2}X=-X,\qquad g^{\prime}(JX,JY)=-g^{\prime}(X,Y)
\end{equation*}
for all vector fields $X,Y \in \mathfrak{X}(M^{\prime})$ (the Lie
algebra of the differentiable vector fields on $M^{\prime}$). The
associated metric $\widetilde{g^{\prime}}$ of the manifold is given
by $\widetilde{g^{\prime}}(X{,}Y)=g^{\prime}(X{,}JY).$ Both metrics
are necessarily of signature $(n^{\prime}{,}n^{\prime}).$

Further, $X,Y,Z,U$ will stand for arbitrary differentiable vector
fields on the manifold, and $x,y,z,u$ -- arbitrary vectors in its
tangent space at an arbitrary point.

The (0,3)-tensor $F^{\prime}$ on $M^{\prime}$ is defined by
$F^{\prime}(X{,}Y{,}Z)=g^{\prime}((\nabla_{X}^{\prime}\,J)Y{,}Z)$,
where $\nabla^{\prime}$ is the Levi-Civita connection of
$g^{\prime}$.

A decomposition to three basic classes of the considered manifolds
with respect to $F ^{\prime}$ is given in \cite{gabo:note}. In this
paper we shall consider only the class $\mathcal{W}_{0} : \
F^{\prime}=0$ of the K\"{a}hler manifolds with Norden metric. The
complex structure $J$ is parallel on every
$\mathcal{W}_{0}$-manifold, i.e. $\nabla^{\prime}J=0$.

The curvature tensor field $R^{\prime}$, defined by
$R^{\prime}(X,Y)Z=\nabla_{X}^{\prime}\nabla_{Y}^{\prime}Z-\nabla_{Y}^{\prime
}\nabla_{X}^{\prime}Z-\nabla_{\lbrack X,Y]}^{\prime}Z$, has the
property $R^{\prime}(X,Y,Z,U)=-R^{\prime}(X,Y,JZ,JU)$ on a
$\mathcal{W}_{0}$-manifold. Using the first Bianchi identity and the
last property of $R$ it
follows $R^{\prime}(X,JY,JZ,U)=-R^{\prime}(X,Y,Z,U).$ Therefore, the tensor field $\widetilde{R^{\prime}}:\widetilde{R^{\prime}}%
(X,Y,Z,U)=R^{\prime}(X,Y,Z,JU)$ has the properties of a K\"{a}hler
curvature tensor and it is called \emph{an associated curvature
tensor}.

The essential curvature-like tensors are defined by:
\begin{equation*}
\begin{split}
\pi_{1}^{\prime}(x,y,z,u) &
=g^{\prime}(y,z)g^{\prime}(x,u)-g^{\prime
}(x,z)g^{\prime}(y,u), \\
\pi_{2}^{\prime}(x,y,z,u) &
=g^{\prime}(y,Jz)g^{\prime}(x,Ju)-g^{\prime
}(x,Jz)g^{\prime}(y,Ju), \\
\pi_{3}^{\prime}(x,y,z,u) &
=-g^{\prime}(y,z)g^{\prime}(x,Ju)+g^{\prime
}(x,z)g^{\prime}(y,Ju) \\
&
{\quad}-g^{\prime}(y,Jz)g^{\prime}(x,u)+g^{\prime}(x,Jz)g^{\prime}(y,u).
\end{split}
\end{equation*}

For every nondegenerate section $\alpha^{\prime}$ in $T_{p^{\prime}}M^{%
\prime }$, $p^{\prime}\in M^{\prime}$, with a basis $\left\{
x,y\right\} $
there are known the following sectional curvatures 
\cite{boga:curv}:
$k^{\prime}(\alpha^{\prime};p^{\prime})=k^{\prime}(x,y)=\frac{R^{\prime
}(x,y,y,x)}{\pi_{1}^{\prime}(x,y,y,x)}$ -- the usual Riemannian
sectional curvature; $\widetilde{k^{\prime}}(\alpha^{\prime};p^{\prime})=\widetilde{k^{\prime}}%
(x,y)=\frac{\widetilde{R^{\prime}}(x,y,y,x)}{\pi_{1}^{\prime}(x,y,y,x)}$
-- an associated sectional curvature.

The sectional curvatures of an arbitrary holomorphic section
$\alpha^{\prime} $ (i.e. $J\alpha^{\prime}=\alpha^{\prime}$) is zero
on a K\"{a}hler manifold with Norden metric \cite{boga:curv}.

For the totally real sections $\alpha^{\prime}$(i.e.
$J\alpha^{\prime}\perp \alpha^{\prime}$) it is proved the following

\begin{theorem}\emph{(\cite{boga:curv})}
\label{t1} Let $M^{\prime}$ ($2n^{\prime}\geq4$) be a K\"{a}hler
manifold with
Norden metric. $M^{\prime}$ is of constant totally real sectional curvatures $%
\nu^{\prime}$and $\widetilde{\nu^{\prime}}$, i.e. $k^{\prime}(\alpha^{%
\prime};p^{\prime})=\nu^{\prime}(p^{\prime})$,\ $\widetilde{k^{\prime}}%
(\alpha^{\prime};p^{\prime})=\widetilde{\nu^{\prime }}(p^{\prime})$
whenever
$\alpha^{\prime}$ is a nondegenerate totally real section in $%
T_{p^{\prime}}M^{\prime}$, $p^{\prime}\in M^{\prime}$, if and only
if
\begin{equation*}
R^{\prime}=\nu^{\prime}\left[ \pi_{1}^{\prime}-\pi_{2}^{\prime}\right] +%
\widetilde{\nu^{\prime}}\pi_{3}^{\prime}.   \label{4}
\end{equation*}
\noindent Both functions $\nu^{\prime}$and
$\widetilde{\nu^{\prime}}$are constant if $M^{\prime}$is connected
and $2n^{\prime}\geq6$.
\end{theorem}

\subsection{Almost contact manifolds with Norden metric}

Let $(M,\varphi,\xi,\eta,g)$ be a $(2n{+}1)$-dimensional almost
contact manifold with Norden metric, i.e. $(\varphi,\xi,\eta)$ is an
almost contact structure determined by a tensor field $\varphi$ of
type $(1,1)$, a vector field $\xi$ and an $1$-form $\eta$ on $M$
satisfying the conditions:
\begin{equation*}
\varphi^{2}X=-X+\eta(X)\xi,\qquad\eta(\xi)=1,   \label{5}
\end{equation*}
and in addition the almost contact manifold $(M,\varphi,\xi,\eta)$
admits a metric $g$ such that \cite{gamigri:alcont}
\begin{equation*}
g(\varphi X,\varphi Y)=-g(X,Y)+\eta(X)\eta(Y).   \label{6}
\end{equation*}

There are valid the following immediate corollaries:
$\eta\circ\varphi=0$, $\varphi\xi=0$, $\eta(X)=g(X{,}\xi)$,
$g(\varphi X{,}Y)=g(X{,}\varphi Y)$.

The associated metric $\widetilde{g}$ given by $\widetilde{g}(X{,}Y)=g(X{,}%
\varphi Y)+\eta(X)\eta(Y)$ is a Norden metric, too. Both metrics are
indefinite of signature $(n,n+1)$.

The Levi-Civita connection of $g$ will be denoted by $\nabla$. The
tensor field $F$ of type (0{,}3) on $M$ is defined by
$F(X{,}Y{,}Z)=g((\nabla_{X}\,\varphi)Y{,}Z).$

If $\{e_{i},\xi\}\,(i=1,2,\dots,2n)$ is a basis of $T_{p}M$ and
$(g^{ij})$ is the inverse matrix of $(g_{ij})$, then the following
$1$-forms are
associated with $F$:%
\begin{equation*}
\theta(\cdot)=g^{ij}F(e_{i},e_{j},\cdot),\quad \theta^{\ast}(\cdot
)=g^{ij}F(e_{i},\varphi
e_{j},\cdot),\quad\omega(\cdot)=F(\xi,\xi,\cdot).
\end{equation*}

A classification of the almost contact manifolds with
Norden metric with respect to $F$ is given in \cite{gamigri:alcont}, where eleven basic classes $\mathcal{F}%
_{i}$ are defined. In the present paper we consider the following
classes:

\begin{equation}
\begin{split}
\mathcal{F}_{4}:\; & F(x,y,z)=-\frac{\theta(\xi)}{2n}\{g(\varphi
x{,}\varphi
y)\eta(z){+}g(\varphi x{,}\varphi z)\eta(y)\}; \\
\mathcal{F}_{5}:\; &
F(x,y,z)=-\frac{\theta^{\ast}(\xi)}{2n}\{g(x,\varphi
y)\eta(z)+g(x,\varphi z)\eta(y)\}; \\
\mathcal{F}_{6}:\; & F(x,y,z)=F(x,y,\xi)\eta(z)+F(x,\xi,z)\eta(y),
\quad\theta(\xi)=\theta^{\ast}(\xi)=0, \\
& F(x,y,\xi)=F(y,x,\xi),\quad F(\varphi x,\varphi y,\xi)=-F(x,y,\xi); \\
\mathcal{F}_{11}:\; & F(x,y,z)=\eta (x)\{\eta (y)\omega (z)+\eta
(z)\omega (y)\}.
\end{split}
\label{9}
\end{equation}

The classes $\mathcal{F}_{i}\oplus\mathcal{F}_{j}$, etc., are
defined in a
natural way by the conditions of the basic classes. The special class $%
\mathcal{F}_{0}$ : $F=0$ is contained in each of the defined
classes. The $\mathcal{F}_{i}^{0}$-manifold is an $\mathcal{F}_{i}$
-manifold ($i=1,4,5,11$) with closed $1$-forms $\theta$, $\theta^*$
and $\omega\circ\varphi$.

The following tensors are essential curvature tensors on $M$:%
\begin{equation*}
\begin{split}
\pi_{1}(x,y,z,u) & =g(y,z)g(x,u)-g(x,z)g(y,u), \\
\pi_{2}(x,y,z,u) & =g(y,\varphi z)g(x,\varphi u)-g(x,\varphi z)g(y,\varphi
u), \\
\pi_{3}(x,y,z,u) & =-g(y,z)g(x,\varphi u)+g(x,z)g(y,\varphi u) \\
& {\quad}-g(y,\varphi z)g(x,u)+g(x,\varphi z)g(y,u), \\
\pi_{4}(x,y,z,u) & =\eta(y)\eta(z)g(x,u)-\eta(x)\eta(z)g(y,u) \\
& {\quad}+\eta(x)\eta(u)g(y,z)-\eta(y)\eta(u)g(x,z), \\
\pi_{5}(x,y,z,u) & =\eta(y)\eta(z)g(x,\varphi u)-\eta(x)\eta(z)g(y,\varphi u)
\\
& {\quad}+\eta(x)\eta(u)g(y,\varphi z)-\eta(y)\eta(u)g(x,\varphi z).
\end{split}%
\end{equation*}

In \cite{mangri:tens} it is established that the tensors $%
\pi_{1}-\pi_{2}-\pi_{4}$ and $\pi_{3}+\pi_{5}$ are K\"{a}hlerian,
i.e. they have the condition of a curvature-like tensor $L$:
$L(X,Y,Z,U)=-L(X,Y,\varphi Z,\varphi U)$.

Let $R$ be the curvature tensor of $\nabla.$ The tensors $R$ and $\widetilde{%
R}:$ $\widetilde{R}(x,y,z,u)=R(x,y,z,\varphi u)$ are K\"{a}hlerian on any $%
\mathcal{F}_{0}$-manifold.

There are known the following sectional curvatures with respect to
$g$ and $R
$ for every nondegenerate section $\alpha$ in $T_{p}M$ with a basis $\{x,y\}$%
:%
\begin{equation*}
k(\alpha;p)=k(x,y)=\frac{R(x,y,y,x)}{\pi_{1}(x,y,y,x)},\quad\widetilde {k}%
(\alpha;p)=\widetilde{k}(x,y)=\frac{\widetilde{R}(x,y,y,x)}{\pi _{1}(x,y,y,x)%
}.
\end{equation*}
In \cite{nakgri:subcotwo} there are introduced the following
special sections in $T_{p}M$: a
$\xi$-section (e.g. $\{\xi,x\}$), a $\varphi$-holomorphic section (i.e. $%
\alpha=\varphi\alpha$) and a totally real section (i.e. $\alpha\perp
\varphi\alpha$).

The canonical curvature tensor $K$ is introduced in
\cite{mangri:tens}. The tensor $K$ is a curvature tensor with
respect to the canonical connection $D$
defined by%
\begin{equation}
D_{X}Y=\nabla_{X}Y+\frac{1}{2}\left\{ (\nabla_{X}\varphi)\varphi
Y+(\nabla_{X}\eta)Y.\xi\right\} -\eta(Y)\nabla_{X}\xi.   \label{12a}
\end{equation}
The connection $D$ is a natural connection, i.e. the structural
tensors are
parallel with respect to $D.$ Let us note that the tensor $K$ out of $%
\mathcal{F}_{0}$ has the properties of $R$ in $\mathcal{F}_{0}.$

\section{Curvatures on
the real time-like hypersurfaces of a K\"ahler manifold with Norden
metric}

In \cite{man:hyp} two types of real hypersurfaces of a complex
manifold with Norden metric are introduced. The obtained
submanifolds are almost contact manifolds with Norden metric. Let
us recall the real time-like hypersurface with respect to the
Norden metric.

The hypersurface $M$ of an almost complex manifold with Norden metric $%
(M^{\prime}, J, g^{\prime})$, determined by the condition the normal
unit $N$ to be time-like regarding $g ^{\prime}$ (i.e.
$g^{\prime}(N,N)=-1$), equipped with the almost contact structure
with Norden metric
\begin{equation}\label{1type}
\begin{array}{rl}
& \varphi :=J+ \cos t .g^{\prime}(\cdot ,JN)\{\cos t .N- \sin t .JN\},\smallskip \\
& \xi := \sin t .N+ \cos t .JN , \quad \eta:= \cos t .
g^{\prime}(\cdot , JN) , \quad g :=g ^{\prime}\vert_M ,
\end{array}
\end{equation}
where $t:=\arctan \left\{ g^{\prime}(N,JN) \right\}$ for $t \in
\left( - \frac{\pi}{2};\frac{\pi}{2}\right) $, is called \emph{a
real time-like hypersurface} of $(M^{\prime}, J, g^{\prime})$.

In the case when $(M^{\prime},J,g^{\prime})$ is a K\"{a}hler
manifold with Norden metric (i.e. a $\mathcal{W}_{0}$-mani\-fold),
in \cite{man:timehyp} it is ascertained the following statement:
The class $\mathcal{F}_{4} \otimes \mathcal{F}_{5} \otimes
\mathcal{F}_{6} \otimes \mathcal{F}_{11}$ is the class of the real
time-like hypersurfaces of a K\"{a}hler manifold with Norden
metric. There are 16 classes of these hypersurfaces in all. When
$n = 1$ the class $\mathcal{F}_ 6$ is restricted to $\mathcal{F}_
0$. Therefore, for a 4-dimensional K\"{a}hler manifold with Norden
metric there are only 8 classes of the considered hypersurfaces.

The tensor $F$ and the second fundamental tensor $A$ of the
considered type of hypersurfaces have the following form,
respectively:
\begin{equation}  \label{23}
\begin{array}{l}
F(X,Y,Z) =\sin t \left\{ g(AX, \varphi Y) \eta (Z) + g(AX, \varphi
Z) \eta (Y) \right\}\smallskip \\- \cos t \left\{ g(AX,Y) \eta (Z) +
g(AX,Z) \eta (Y)-2 \eta (AX) \eta (Y) \eta (Z) \right\},\smallskip\\
AX= - \frac{dt(\xi )}{2 \cos t} \eta (X) \xi  -\sin t \{ \nabla_X
\xi +
g(\nabla_\xi \xi ,X) \xi \}\smallskip \\
\qquad +\cos t \{ \varphi\nabla_X \xi + g(\varphi\nabla_\xi \xi ,X)
\xi \}.
\end{array}
\end{equation}
The basic classes of the considered hypersurfaces are characterized
in terms of the second fundamental tensor by the conditions
\cite{man:hyp}:
\begin{equation*}
\begin{array}{rl}
\mathcal{F}_0 \, : \;\; & A=- \frac{dt( \xi )}{2 \cos t} \eta\otimes\xi ;\smallskip \\
\mathcal{F}_4 \, : \;\; & A=- \frac{dt( \xi )}{2 \cos t} \eta\otimes\xi - \frac{%
\theta (\xi )}{2n} \{ \sin t. \varphi - \cos t. \varphi^2 \} ; \\[4pt]
\mathcal{F}_5 \, : \;\; & A=- \frac{dt( \xi )}{2 \cos t} \eta\otimes\xi + \frac{%
\theta^* (\xi )}{2n} \{ \cos t. \varphi + \sin t. \varphi^2 \} ; \\[4pt]
\mathcal{F}_6 \, : \;\; & A \circ\varphi = \varphi\circ A \, , \;\;
tr A - \frac{dt(
\xi )}{2 \cos t}=tr (A \circ\varphi )=0 ; \\[4pt]
\mathcal {F}_{11} \, : \;\; & A=- \frac{dt( \xi )}{2 \cos t}
\eta\otimes\xi - \cos t
\{\eta\otimes\Omega + \omega\otimes\xi\} \\[4pt]
& \qquad\qquad\quad\qquad\hspace{0.1in} - \sin t
\{\eta\otimes\varphi\Omega + (\omega\circ\varphi)
\otimes\xi\},\hspace{0.1in} \omega(\cdot) = g(\cdot,\Omega).

\end{array}
\end{equation*}

According to the formulas of Gauss and Weingarten in this case
$\nabla^{\prime}_X Y=\nabla_X Y-g(AX,Y)N$, $\nabla^{\prime}_X
N=-AX$, we get the relation between the curvature tensors
$R^{\prime}$ and $R$ of the $\mathcal{W}_{0}$-manifold
$(M^{\prime},J,g^{\prime})$ and its hypersurface
$(M,\varphi,\xi,\eta,g)$, respectively:
\begin{equation*}
\begin{array}{l}
R^{\prime}(x,y,z,u)=R(x,y,z,u)+\pi_{1}(Ax,Ay,z,u),\smallskip\\R^{\prime}(x,y)N=
-\left(\nabla_{x}A\right)y+\left(\nabla_{y}A\right) x.
\end{array}
\end{equation*}

Hence, having in mind Theorem \ref{t1}, we obtain:%
\begin{equation*}
\begin{array}{l}
R(x,y,z,u) =\left\{ \nu^{\prime}\left[
\pi_{1}^{\prime}-\pi_{2}^{\prime } \right]
+\widetilde{\nu^{\prime}}\pi_{3}^{\prime}\right\} (x,y,z,u)-\pi
_{1}(Ax,Ay,z,u),\smallskip \\
R(x,y,\varphi z,\varphi u) = -\left \{R - \nu^{\prime}[\pi_{4}-\tan
t \pi_{5}]+ \widetilde{\nu^{\prime}}[\pi_{5}+\tan t
\pi_{4}]\right\}(x,y,z,u)\smallskip\\
\qquad\qquad\qquad\hspace{0.1in}
-[\pi_{1}+\pi_{2}](Ax,Ay,z,u),\smallskip\\
R(x,y)\xi =\left\{\nu^{\prime}[\pi_{4}-\tan t \pi_{5}] -
\widetilde{\nu^{\prime}}[\pi_{5}+\tan t \pi_{4}]\right\}(x,y)\xi -
\pi_{1}(Ax,Ay)\xi,\smallskip\\
R(x,y)N = -\frac{1}{\cos t}[\nu^{\prime}\pi_{5} +
\widetilde{\nu^{\prime}}\pi_{4}](x,y)\xi.
\end{array}
\end{equation*}
Therefore
\begin{equation}
(\nabla_{x}A)y - (\nabla_{y}A)x = \frac{1}{\cos t}[\nu^{\prime}
\pi_{5} + \widetilde{\nu^{\prime}} \pi_{4}](x,y)\xi. \label{14}
\end{equation}

Having in mind the equations: $g^{\prime}(y,Jz)=g(y,\varphi z)+\tan
t\ \eta(y)\eta(z) $, $\pi_{1}^{\prime}=\pi_{1}$, $\pi
_{2}^{\prime}=\pi_{2}+ \tan t\ \pi_{5}$, $\pi_{3}^{\prime}=\pi_{3}-
\tan t\ \pi_{4}$, which are valid for real time-like hypersurfaces,
we obtain
\begin{proposition}
\label{l1} A real time-like hypersurface of a K\"{a}hler manifold
with Norden metric of constant totally real sectional curvatures
$\nu^{\prime}$ and $\widetilde{\nu^{\prime}}$ has the following
curvature properties:
\begin{equation*}
\begin{array}{l}
R(x,y,z,u)=\left\{ \nu^{\prime}[\pi_{1}-\pi_{2}-\tan t \pi_{5}]+
\widetilde{\nu^{\prime}}[\pi_{3}-\tan t \pi_{4}]\right\}(x,y,z,u) \smallskip\\
\qquad\qquad\hspace{0.18in} - \pi_{1}(Ax,Ay,z,u),\smallskip\\
\tau = 4n^{2}\nu^{\prime} - 4n\widetilde{\nu^{\prime}}\tan t - (\tr
A)^{2} + \tr A^{2},\smallskip\\
\widetilde{\tau} = -2n\nu^{\prime}\tan t +
2n(2n-1)\widetilde{\nu^{\prime}} - \tr A \tr (A\circ \varphi) +
\tr(A^{2}\circ \varphi);
\end{array}
\end{equation*}
for a $\xi$-section $\left\{ \xi,x\right\} $
\begin{equation*}
k(\xi,x)=\nu^{\prime} -\widetilde{\nu^{\prime}}\tan t
-[\nu^{\prime}\tan t + \widetilde{\nu^{\prime}}]\frac{g(x,\varphi
x)}{g(x,x)-\eta(x)^{2}} -
\frac{\pi_{1}(A\xi,Ax,x,\xi)}{g(x,x)-\eta(x)^{2}};
\end{equation*}
for a $\varphi$-holomorphic section $\left\{ \varphi
x,\varphi^{2}x\right\} $ and for a totally real section $\left\{
x,y\right\} $, orthogonal to $\xi $, respectively:
\begin{equation*}
k(\varphi x,\varphi^{2}x)=-\frac{\pi_{1}(A\varphi x,
A\varphi^{2}x,\varphi^{2}x,\varphi x)}{\pi_{1}(\varphi x,
\varphi^{2}x,\varphi^{2}x,\varphi x)},\quad k(x,y)=\nu^{\prime} -
\frac{\pi_{1}(Ax,Ay,y,x)}{\pi_{1}(x,y,y,x)}.
\end{equation*}
\end{proposition}

If $(M,\varphi,\xi,\eta,g)$ is a real time-like hypersurface of
$\mathcal{W}_{0}$-manifold,  then (\ref{12a}), (\ref{23}) and
(\ref{14}) imply that the canonical curvature tensor has the form
\begin{equation*}
\begin{array}{c}
K(x,y,z,u)=R(x,y,\varphi^{2}z,\varphi^{2}u)+\pi_{1}(Ax,Ay,\varphi
z,\varphi u)\smallskip\\
+\sin t \left\{ \sin t [\pi_{1}-\pi_{2}-\pi_{4}] - \cos t
[\pi_{3}+\pi_{5}] \right\}(Ax,Ay,z,u).
\end{array}
\end{equation*}


Then, because of the last equation and Proposition \ref{l1} we have
\begin{proposition}
Let $(M,\varphi,\xi,\eta,g)$ be a real time-like hypersurface of a $\mathcal{%
W}_{0}$-manifold $(M^{\prime},J,g^{\prime})$ of constant totally
real sectional curvatures. Then $K$ of $M$ 
is K\"{a}hlerian and
\begin{equation*}
\begin{array}{l}
K(x,y,z,u) =\left\{ \nu^{\prime}[\pi_{1}-\pi_{2}-\pi_{4}]+
\widetilde{\nu^{\prime}}[\pi_{3} + \pi_{5}]\right\}(x,y,z,u)\smallskip\\
\qquad\qquad\hspace{0.21in} - \cos t\big\{ \cos t
[\pi_{1}-\pi_{2}-\pi_{4}]+\sin t[\pi_{3} +
\pi_{5}]\big\}(Ax,Ay,z,u),\smallskip\\
\tau(K) = 4n(n-1)\nu^{\prime} - \cos t (a\cos t + 2b\sin
t),\smallskip\\
\widetilde{\tau}(K) = 4n(n-1)\widetilde{\nu^{\prime}} - \cos t (a
\sin t - 2b\cos t),\medskip\\
a = (\tr A)^{2} - \tr A^{2} - [\tr (A \circ \varphi)]^{2}+ \tr (A
\circ \varphi)^{2} - 2 \eta(A\xi) \tr A +
2g(A\xi,A\xi),\smallskip\\
b= \tr(A^{2}\circ \varphi) - \tr A \tr(A\circ \varphi) +
\eta(A\xi)\tr(A\circ\varphi) - g(\varphi A\xi,A\xi).
\end{array}
\end{equation*}
\end{proposition}

\section{Curvatures on $\mathcal{W}_{0}$'s real time-like hypersurfaces,
belonging to the main classes}

Now, let $(M,\varphi,\xi,\eta,g)$ belong to the widest integrable main class $\mathcal{F%
}_{4}\oplus\mathcal{F}_{5}$ of the real time-like hypersurfaces. Let
us recall that a class of almost contact manifolds with Norden
metric is said to be main if the tensor $F$ is expressed explicitly by the structural tensors $%
\varphi,\xi,\eta,g$. In this case for the second fundamental tensor
we have\cite{man:hyp}:
\begin{equation}\label{14a}
\begin{array}{l}
A=-\frac{dt(\xi)}{2\cos t}\eta \otimes \xi -\frac{1}{2n}\big\{
[\theta(\xi)\sin t - \theta^{\ast}(\xi)\cos t]\varphi \smallskip\\
\qquad - \left[\theta(\xi)\cos t + \theta^{\ast}(\xi)\sin t
\right]\varphi^{2}\big\},\smallskip\\
\tr A = -\frac{dt(\xi)}{2\cos t} - \theta(\xi)\cos t -
\theta^{\ast}(\xi)\sin t,\smallskip\\ \tr (A\circ \varphi) =
\theta(\xi)\sin t - \theta^{\ast}(\xi)\cos t.
\end{array}
\end{equation}
Then, having in mind the last identities and Proposition \ref{l1},
we obtain
\begin{corollary}
\label{c1} If a real time-like hypersurface of a K\"{a}hler
manifold with Norden metric of constant totally real sectional
curvatures is an $\mathcal{F}_{4}\oplus\mathcal{F}_{5}$-manifold,
then it has the following curvature properties:
\begin{equation*}
\begin{array}{l}
R= \nu^{\prime}[\pi_{1} - \pi_{2} - \tan t \pi_{5}] +
\widetilde{\nu^{\prime}}[\pi_{3} - \tan t \pi_{4}]\smallskip\\
\quad - \frac{dt(\xi)}{4n\cos t}\big\{ \theta(\xi)[\sin t \pi_{5} +
\cos t \pi_{4}] + \theta^{\ast}(\xi)[\sin t \pi_{4} - \cos t
\pi_{5}]\big\}\smallskip\\ \quad - \frac{\theta(\xi)^{2} +
\theta^{\ast}(\xi)^{2}}{4n^{2}}\pi_{2} -\frac{\big(\theta(\xi)\cos t
+ \theta^{\ast}(\xi)\sin
t\big)^{2}}{4n^{2}}[\pi_{1}-\pi_{2}-\pi_{4}]\smallskip\\
\quad+\frac{\big(\theta(\xi)\cos t + \theta^{\ast}(\xi)\sin
t\big)\big(\theta(\xi)\sin t - \theta^{\ast}(\xi)\cos
t\big)}{4n^{2}}[\pi_{3}+\pi_{5}],\smallskip\\
\tau = 4n(n\nu^{\prime} -\widetilde{\nu^{\prime}}\tan t) -
dt(\xi)\theta(\xi) - dt(\xi)\theta^{\ast}(\xi)\tan t\smallskip\\
\quad - \frac{n-1}{n}(\theta(\xi)\cos t + \theta^{\ast}(\xi)\sin
t)^{2} -
\frac{\theta(\xi)^{2}+\theta^{\ast}(\xi)^{2}}{2n},\smallskip\\
\widetilde{\tau}=2n(n-1)\widetilde{\nu^{\prime}} +
2n\nu^{\prime}\tan t + \frac{dt(\xi)\theta(\xi)}{2}\tan t -
\frac{dt(\xi)\theta^{\ast}(\xi)}{2} \smallskip\\
\quad + \frac{n-1}{n}\big(\theta(\xi)\sin t - \theta^{\ast}(\xi)\cos
t)(\theta(\xi)\cos t + \theta^{\ast}(\xi)\sin t \big),\smallskip\\
k(\varphi
x,\varphi^{2}x)=-\frac{\theta(\xi)^{2}+\theta^{\ast}(\xi)^{2}}{4n^{2}},\quad
k(x,y)=\nu^{\prime} - \frac{\big(\theta(\xi)\cos t +
\theta^{\ast}(\xi)\sin t\big)^{2}}{4n^{2}}.
\end{array}
\end{equation*}
\end{corollary}
Let us remark that we can obtain the corresponding properties for
the classes $\mathcal{F}_{4}$,
$\mathcal{F}_{5}$ and $\mathcal{F}_{0}$, if we substitute $\theta^{\ast}%
(\xi)=0$, $\theta(\xi)=0$ and $\theta(\xi)=\theta^{\ast}(\xi)=0$, respectively.

Using the equations (\ref{12a}) and (\ref{9}), we express the
canonical connection explicitly for the class
$\mathcal{F}_{4}\oplus\mathcal{F}_{5}$ as follows
\begin{equation*}
D_{X}Y=\nabla_{X}Y+\frac{\theta(\xi)}{2n}\left\{ g(x,\varphi y)\xi
-\eta(y)\varphi x\right\} -\frac{\theta^{\ast}(\xi)}{2n}\left\{
g(\varphi x,\varphi y)\xi-\eta(y)\varphi^{2}x\right\}.
\end{equation*}
Let $(M,\varphi,\xi ,\eta,g)\in%
\mathcal{F}_{4}^{0}\oplus\mathcal{F}_{5}^{0},$ i.e. $(M,\varphi
,\xi,\eta,g)$ is an $\left(
\mathcal{F}_{4}\oplus\mathcal{F}_{5}\right) $-manifold with closed
$1$-forms $\theta$ and $\theta^{\ast}.$ The canonical curvature
tensor $K$ of any $\left( \mathcal{F}_{4}^{0}\oplus \mathcal{F}%
_{5}^{0}\right) $-manifold is K\"{a}hlerian and it has the form
\begin{equation*}
\begin{array}{l}
K=R+\frac{\xi\theta(\xi)}{2n}\pi_{5}+\frac{\xi\theta^{\ast}(\xi)}{2n}
\pi_{4}\smallskip \\
\qquad+\frac{\theta^{2}(\xi)}{4n^{2}}\left[ \pi_{2}-\pi_{4}\right] +\frac{%
\theta^{\ast2}(\xi)}{4n^{2}}\pi_{1}-\frac{\theta(\xi)\theta^{\ast}(\xi)}{%
4n^{2}}\left[ \pi_{3}-\pi_{5}\right].
\end{array}
\end{equation*}
Then, using Corollary \ref{c1}, we ascertain the truthfulness of the
following
\begin{corollary}
\label{c2} If a real time-like hypersurface of a K\"{a}hler
manifold with Norden metric of constant totally real sectional
curvatures is an
$\bigl(\mathcal{F}_{4}^{0}\oplus\mathcal{F}_{5}^{0}\bigr)$-manifold,
then $K$ is expressed in the following way:
\begin{equation*}
\begin{array}{l}
\begin{array}{l}
K=\big(\nu^{\prime} +
\frac{\theta^{\ast}(\xi)}{4n^{2}}\big)[\pi_{1}-\pi_{2}] +
\big(\widetilde{\nu^{\prime}}-\frac{\theta(\xi)\theta^{\ast}(\xi)}{4n^{2}}\big)\pi_{3}\smallskip\\
\quad - \big(\widetilde{\nu^{\prime}}\tan
t+\frac{dt(\xi)\theta(\xi)}{4n}+\frac{dt(\xi)\theta^{\ast}(\xi)}{4n}\tan
t - \frac{\xi
\theta^{\ast}(\xi)}{2n}+\frac{\theta(\xi)^{2}}{4n^{2}}\big)\pi_{4}\smallskip\\
\quad - \big(\nu^{\prime}\tan t + \frac{dt(\xi)\theta(\xi)}{4n}\tan
t - \frac{dt(\xi)\theta^{\ast}(\xi)}{4n}-\frac{\xi
\theta(\xi)}{2n}-\frac{\theta(\xi)\theta^{\ast}(\xi)}{4n^{2}}\big)\pi_{5}\smallskip\\
\quad - \frac{\bigl(\theta(\xi)\cos t + \theta^{\ast}(\xi)\sin
t\bigr)^{2}}{4n^{2}}[\pi_{1} - \pi_{2}-\pi_{4}]\smallskip\\
\quad + \frac{\bigl(\theta(\xi)\sin t - \theta^{\ast}(\xi)\cos
t\bigr)\bigl(\theta(\xi)\cos t + \theta^{\ast}(\xi)\sin
t\bigr)}{4n^{2}}[\pi_{3}+\pi_{5}].
\end{array}
\end{array}
\end{equation*}
\end{corollary}

We compute the expression $\left( \nabla_{x}A\right) y-\left( \nabla
_{y}A\right) x$ using (\ref{14a}) and we compare the result with
(\ref{14}). Thus, we get the relations
\begin{equation}
\begin{array}{l}
\nu^{\prime} = -\frac{dt(\xi)\theta(\xi)}{4n} + \frac{\cos
t}{2n}\big[\xi \theta(\xi)\sin t - \xi \theta^{\ast}(\xi)\cos
t\big]\smallskip\\ \qquad + \frac{\cos
^{2}t}{4n^{2}}\big[\theta(\xi)^{2} - \theta^{\ast}(\xi)^{2}\big]+
\frac{\sin t\cos t}{2n^{2}}\theta(\xi)\theta^{\ast}(\xi)\smallskip\\
\qquad  + \frac{dt(\xi)}{2n}\cos t
\big[\theta(\xi)\cos t + \theta^{\ast}(\xi)\sin t\big],\smallskip\\
\widetilde{\nu^{\prime}} = -\frac{dt(\xi)\theta^{\ast}(\xi)}{4n} +
\frac{\cos t}{2n}\big[\xi \theta(\xi) \cos t + \xi
\theta^{\ast}(\xi) \sin t\big]\smallskip\\ \qquad + \frac{\sin t
\cos t}{4n^{2}}\big[\theta^{\ast}(\xi)^{2} - \theta(\xi)^{2}\big]+
\frac{\cos ^{2} t}{2n^{2}}\theta(\xi)\theta^{\ast}(\xi)\smallskip\\
\qquad  - \frac{dt(\xi)}{2n}\cos t \big[\theta(\xi)\sin t -
\theta^{\ast}(\xi)\cos t\big].
\end{array}  \label{15a}
\end{equation}
Hence, we have
\begin{equation*}
\begin{array}{l}
K = \lambda [\pi_{1} - \pi_{2} - \pi_{4}] + \mu [\pi_{3} +
\pi_{5}],\smallskip\\
R = \lambda [\pi_{1} - \pi_{2} - \pi_{4}] +
\mu [\pi_{3} + \pi_{5}]
-\frac{\xi \theta^{\ast}(\xi)}{2n}\pi_{4} - \frac{\xi \theta(\xi)}{2n}\pi_{5} \smallskip\\
\phantom{R = }
-\frac{\theta^{\ast}(\xi)^{2}}{4n^{2}}\pi_{1}-\frac{\theta(\xi)^{2}}{4n^{2}}[\pi_{2}
- \pi_{4}] +
\frac{\theta(\xi)\theta^{\ast}(\xi)}{4n^{2}}[\pi_{3}-\pi_{5}],\medskip\\
\lambda = -\frac{dt(\xi)\theta(\xi)}{4n} + \frac{dt(\xi)}{2n}\cos t
\big[\theta(\xi)\cos t +
\theta^{\ast}(\xi)\sin t\big] \smallskip\\
\phantom{\lambda =}
 + \frac{\cos t}{2n}\big[\xi \theta(\xi) \sin t -
\xi
\theta^{\ast}(\xi)\cos t\big],\smallskip\\
\mu = -\frac{dt(\xi)\theta^{\ast}(\xi)}{4n}- \frac{dt(\xi)}{2n}\cos
t \big[\theta(\xi)\sin t - \theta^{\ast}(\xi)\cos t\big]
\smallskip\\
\phantom{\mu =}
+ \frac{\cos t}{2n}\big[\xi \theta(\xi) \cos t + \xi
\theta^{\ast}(\xi)\sin t\big].
\end{array}
\end{equation*}

We solve the system (\ref{15a}) with respect to the functions
$\theta(\xi)$ and $\theta^{\ast}(\xi)$ for $t=\textrm{const}$ and
get
\begin{equation*}
\begin{array}{l}
\theta(\xi)=2\varepsilon n\sqrt{\frac{\nu^{\prime}\cos t -
\widetilde{\nu}^{\prime}\sin t + \sqrt{\nu^{\prime
2}+\widetilde{\nu}^{\prime 2}}}{2\cos t}},
\theta^{\ast}(\xi)=2\varepsilon n\frac{\sqrt{\cos
t}(\nu^{\prime}\tan t
+\widetilde{\nu}^{\prime})}{\sqrt{2(\nu^{\prime}\cos t -
\widetilde{\nu}^{\prime}\sin t + \sqrt{\nu^{\prime
2}+\widetilde{\nu}^{\prime 2}})}},
\end{array}
\end{equation*}
where $\varepsilon = \pm 1$. Since $\nu^{\prime}$ and
$\widetilde{\nu^{\prime}}$ are pointwise constant for
$M^{\prime4}$ $(n=1)$ and they are absolute constants for
$M^{\prime2n+2}
$ $(n\geq2)$ (Theorem \ref{t1}), then the functions $\theta(\xi)$ and $%
\theta^{\ast}(\xi),$ which determine the real time-like
hypersurface as an almost contact manifold with Norden metric, are
also pointwise constant on $M^{3}$ and absolute constants on
$M^{2n+1}$ ($n\geq2$). Hence, we have
\begin{theorem}
Let $(M^{\prime},J,g^{\prime})$ be a K\"{a}hler manifold with
Norden metric of constant totally real sectional curvatures. Let
the $\left(
\mathcal{F}_{4}^{0}\oplus\mathcal{F}_{5}^{0}\right)$-manifold
$(M,\varphi,\xi,\eta,g)$, $\dim M \geq 5$, be its real time-like
hypersurface, defined by (\ref{1type}). If $t=\const$, then $K=0$
on $M$ and
\begin{equation*}
\begin{array}{l}
R =-\frac{\theta(\xi)^{2}}{4n^{2}}\left[ \pi_{2}-\pi_{4}\right] -\frac{%
\theta^{\ast}(\xi)^{2}}{4n^{2}}\pi_{1}+\frac{\theta(\xi)\theta^{\ast}(\xi)}{%
4n^{2}}\left[ \pi_{3}-\pi_{5}\right],\smallskip \\
\tau=\frac{\theta(\xi)^{2}}{2n}-\left( 2n+1\right) \frac{
\theta^{\ast}(\xi)^{2}}{2n},\qquad
\widetilde{\tau}=\frac{\theta(\xi)\theta^{\ast}(\xi)}{2n},\smallskip\\
k\left( \xi,x\right)
=\frac{\theta(\xi)^{2}-\theta^{\ast}(\xi)^{2}}{4n^{2}}+
\frac{2\theta(\xi)\theta^{\ast}(\xi)}{4n^{2}}\frac{g(x,\varphi
x)}{g(\varphi x,\varphi x)},\smallskip\\
k\left( \varphi x,\varphi^{2}x\right) =-\frac{\theta(\xi
)^{2}+\theta^{\ast}(\xi)^{2}}{4n^{2}},\qquad k\left( x,y\right)
=-\frac{\theta^{\ast}(\xi)^{2}}{4n^{2}}.
\end{array}
\end{equation*}
\end{theorem}


\end{document}